\documentclass[11pt]{article}
\usepackage[margin=1in]{geometry}
\usepackage[T1]{fontenc}
\usepackage{lmodern}
\usepackage{amsmath,amssymb,amsthm,mathtools}
\usepackage{booktabs,array,enumitem,microtype,tabularx}
\usepackage[numbers,sort&compress]{natbib}
\usepackage{xcolor}
\usepackage{hyperref}
\hypersetup{hidelinks,pdftitle={Matrix--Vector Complexity of Low-Rank Approximation}}

\newcommand{\R}{\mathbb{R}}
\newcommand{\E}{\mathbb{E}}
\newcommand{\Prb}{\mathbb{P}}
\newcommand{\eps}{\varepsilon}
\newcommand{\op}{\mathrm{op}}
\newcommand{\rank}{\operatorname{rank}}
\newcommand{\range}{\operatorname{range}}
\newcommand{\tr}{\operatorname{tr}}
\newcommand{\wtO}{\widetilde O}

\newcommand{\wtTheta}{\widetilde\Theta}
\newcommand{\cF}{\mathcal F}
\newcommand{\cG}{\mathcal G}

\newcommand{\up}{\mathord{\uparrow}}

\newcommand{\OPT}{\operatorname{OPT}}

\newtheorem{theorem}{Theorem}[section]
\newtheorem{lemma}[theorem]{Lemma}
\newtheorem{proposition}[theorem]{Proposition}
\newtheorem{corollary}[theorem]{Corollary}
\theoremstyle{definition}

\newtheorem{remark}[theorem]{Remark}

\usepackage[nameinlink,noabbrev]{cleveref}

\newcolumntype{Y}{>{\raggedright\arraybackslash}X}
\allowdisplaybreaks[1]
\title{\textbf{Matrix--Vector Complexity of Low-Rank Approximation}}
\newcommand{\authornames}{Haihan Zhang, Wendao Wu, Chenheng Zhang, Yanyi Li, Chunyuan Zheng, Cong Fang, Haoxuan Li, Zhouchen Lin}
\newcommand{\affiliation}{Peking University}
\newcommand{\contactemail}{colombo.2026@outlook.com}
\newcommand{\mail}[1]{\href{mailto:#1}{\texttt{#1}}}
\newif\ifanonymous
\anonymousfalse
\ifanonymous
\author{Anonymous Authors}
\hypersetup{pdfauthor={Anonymous Authors}}
\else
\author{Haihan Zhang$^{1,*}$\quad Wendao Wu$^{1,*}$\quad Chenheng Zhang$^{1,*}$\\[0.20em]
Yanyi Li$^1$\quad Chunyuan Zheng$^1$\\[0.20em]
Cong Fang$^{1,\dagger}$\quad Haoxuan Li$^{1,\dagger}$\quad Zhouchen Lin$^{1,\dagger}$\\[0.60em]
\small $^1$\affiliation\\[0.35em]
\footnotesize\mail{\contactemail}\\[0.35em]
\small $^*$Equal contribution.\quad $^\dagger$Corresponding authors.}
\hypersetup{pdfauthor={\authornames}}
\fi
\date{}
\newcommand{\ResearchAgentSystem}{our laboratory's internal auto-research system}

\begin{document}
\maketitle
\begin{abstract}
We establish matching polynomial query bounds for low-rank approximation from exact matrix--vector products. Given an unknown matrix $A\in\R^{m\times n}$, at each step a randomized algorithm chooses either $v\in\R^n$ and receives $Av$, or $u\in\R^m$ and receives $A^\top u$. The choice may depend measurably on all previous queries and replies and on the algorithm's private randomness; each vector product costs one query. The output is a rank-$k$ right projector with Schatten-$p$ residual at most $1+\eps$ times optimal. Write $N=\min\{m,n\}$ and let $Q_p^*$ denote the worst-case query budget for success probability $2/3$ on every input. For every $1\le k<N$ and sufficiently small $\eps$, our lower bounds, combined with existing Krylov upper bounds, give
\[
\begin{aligned}
 Q_p^*&=\wtTheta\!\left(\min\{N,k\min\{p^{1/6}\eps^{-1/3},\eps^{-1/2}\}\}\right)
 &&(2\le p<\infty),\\
 Q_\infty^*&=\wtTheta\!\left(\min\{N,k\eps^{-1/2}\}\right).
\end{aligned}
\]
These bounds have universal constants and allow $p$ to vary with the problem parameters, identifying the transition at $p\eps\asymp1$. A complementary result for each fixed $1\le p<2$ gives $\wtTheta_p(\min\{N,k\eps^{-1/3}\})$, with constants and an accuracy threshold that may depend on $p$. Together, the results recover this fixed-norm rate for every fixed finite $p$, supplying the multiplicative rank dependence missing from previous lower bounds. Tildes suppress logarithmic factors. The proof extends adaptive Wishart deferred decisions to a rectangular factor with a $k$-dimensional nullspace. Posterior overlap gives a short fixed-norm argument, while persistence of small compression eigenvalues controls growing $p$ and the spectral endpoint. Exact range recovery handles target costs of order $k$; the Wishart family covers the remaining regimes.
\end{abstract}

\begin{samepage}
\noindent\textbf{AI Usage.}
Nearly the entire research pipeline for this paper was carried out by
\ResearchAgentSystem{}, powered by GPT-5.6 Sol. The system also conducted a
Lean-backed article audit of the resulting manuscript. The authors subsequently reviewed and approved the
mathematical claims, presentation, and formal artifacts, and take
responsibility for the final manuscript. The complete Lean audit report and the system's technical report will be made public at a later date.
\par
\end{samepage}

\paragraph{Publication Intent and Recognition of Independent Work.}

This preprint is released as a research output of our study of
automated scientific research, for public scrutiny and discussion.
As we do not intend to submit this manuscript to a journal or
conference, it has not undergone the full editorial refinement
and internal review we would normally require for submission.
The authors remain responsible for its contents and for correcting
any errors or inaccurate attribution.

We recognize that releasing these results may affect researchers
who have independently been pursuing the same questions, including
work that is not yet public.
If this applies to your work, we invite you to contact us at
\texttt{colombo.2026@outlook.com}.
A brief description of the overlap, together with any dated
materials you are comfortable sharing, would help us understand
its scope and chronology.
Such materials might include draft excerpts, presentation slides,
or revision records; a complete manuscript or unpublished proof
is not required to begin the conversation.

We will review these cases carefully and update citations,
acknowledgments, and accounts of independent or concurrent work
as appropriate, including revising priority claims where warranted.
Nonpublic materials will be treated confidentially and used only
to assess the overlap unless their authors authorize otherwise.
We will obtain explicit permission from the authors before
publicly discussing such materials or providing them to an
automated research system.

\section{Introduction}
\label{sec:intro}
Given an unknown matrix $A\in\R^{m\times n}$, a target rank $1\le k<\min\{m,n\}$, a Schatten parameter $p\in[1,\infty]$, and an accuracy $\eps>0$, the low-rank approximation problem studied here is to find a rank-$k$ orthogonal projector $P\in\R^{n\times n}$ satisfying
\[
 \|A(I-P)\|_{S_p}\le(1+\eps)\min_{\rank(P_*)=k}\|A(I-P_*)\|_{S_p}.
\]
The minimum is over rank-$k$ orthogonal projectors, and $\|\cdot\|_{S_p}$ denotes the Schatten norm, including the spectral norm at $p=\infty$. We require success probability at least $2/3$ on every input. The algorithm accesses $A$ only through exact products $Av$ or $A^\top u$, chosen adaptively from previous replies. We ask for the smallest worst-case number of such vector products; computation between queries is uncharged.

The principal upper bounds are already known. Bakshi, Clarkson, and Woodruff obtain $\wtO_p(k\eps^{-1/3})$ queries for fixed finite $p$, with explicit norm dependence $\wtO(kp^{1/6}\eps^{-1/3})$ in the finite-$p$ regime \citep{BCW22}. Spectral block Krylov iteration uses $\wtO(k\eps^{-1/2})$ products \citep{MM15}. The fixed-$p$ lower bound established by BCW controls the accuracy dependence at constant rank. For Frobenius approximation, the comparison in Amsel et al.\ explicitly records the additive frontier $\Omega(k+\eps^{-1/3})$ in sufficiently large dimension \citep{AmselEtAl26}. Thus, even in this basic case, the recorded lower bound does not explain the product of rank and accuracy costs.

Our main contribution is a lower bound that recovers this product and extends it through the transition to spectral approximation. Let $N=\min\{m,n\}$. Our main theorem allows $p$ to vary: uniformly for $2\le p<\infty$, the polynomial rate is
\begin{equation}
 \min\{N,k\Phi(p,\eps)\},\qquad
 \Phi(p,\eps)=\min\{p^{1/6}\eps^{-1/3},\eps^{-1/2}\}.
 \label{eq:intro-rate}
\end{equation}
Consequently, the norm dependence grows as $p^{1/6}$ until $p\eps$ is of order one and then reaches the spectral polynomial rate. A complementary bound treats each fixed $1\le p<2$. Combining the two ranges gives the rate $\wtTheta_p(\min\{N,k\eps^{-1/3}\})$ for every fixed finite $p$; this fixed-norm conclusion follows from the two complementary parameter regimes. This conclusion identifies powers of the parameters, rather than exact logarithmic factors. In particular, the dimension-dependent logarithm in the rank-one spectral lower bound of Bakshi and Narayanan remains a finer statement in its regime \citep{BN23}.

\begin{table}[htbp]
\centering
\small
\caption{Matrix--vector query bounds for relative-error right-subspace approximation. Prior lower bounds use their source constant-success conventions and high-dimensional regimes; new bounds use success probability $2/3$, hold for $1\le k<N$, and display rates before dimension saturation: each rate $f$ in the last three rows means $\min\{N,f\}$. Tildes suppress logarithmic factors. Fixed-$p$ constants may depend on $p$; uniform constants are universal. The upper bounds in the last three rows are inherited.}
\label{tab:comparison}
\setlength{\tabcolsep}{4pt}
\begin{tabularx}{\textwidth}{@{}>{\raggedright\arraybackslash}p{0.28\textwidth}YY@{}}
\toprule
Source and regime & Upper bound & Lower bound\\
\midrule
BCW \citep{BCW22}, fixed $p$, $k=1$
 & $\wtO_p(\eps^{-1/3})$ & $\Omega_p(\eps^{-1/3})$\\[5pt]
Recorded frontier \citep{AmselEtAl26}, $p=2$
 & $\wtO(k\eps^{-1/3})$ & $\Omega(k+\eps^{-1/3})$\\[5pt]
MM and BN \citep{MM15,BN23}, $p=\infty$, $k=1$
 & $\wtO(\eps^{-1/2})$ & $\Omega(\eps^{-1/2}\log N)$\\[5pt]
\textbf{This Paper} (\Cref{thm:main})\newline $2\le p<\infty$
 & $\wtO(k\min\{p^{1/6}\eps^{-1/3},\eps^{-1/2}\})$
 & $\boldsymbol{\Omega(k\min\{p^{1/6}\eps^{-1/3},\eps^{-1/2}\})}$\\[5pt]
\textbf{This Paper} (\Cref{thm:main})\newline $p=\infty$
 & $\wtO(k\eps^{-1/2})$
 & $\boldsymbol{\Omega(k\eps^{-1/2})}$\\[5pt]
\textbf{This Paper}\newline (\Cref{prop:small-p}), fixed $1\le p<2$
 & $\wtO_p(k\eps^{-1/3})$
 & $\boldsymbol{\Omega_p(k\eps^{-1/3})}$\\
\bottomrule
\end{tabularx}
\end{table}

The information obstruction is a hidden $k$-dimensional nullspace. We take $W=XX^\top$ with $X\in\R^{d\times(d-k)}$ Gaussian and use the shifted matrix $H=I-W/L$. The shift itself is inherited from the Wishart lower-bound approach of BCW \citep{BCW22}. What changes is the rank deficiency and the information that must be recovered: relative approximation forces the output to capture a constant fraction of an entire nullspace. We extend the adaptive Wishart decomposition of Braverman, Hazan, Simchowitz, and Woodworth \citep{BHSW20} to this rectangular setting and prove
\[
 \E\tr(P\Pi_{\ker W})\le\frac{2k(q+k)}{d-q}
\]
for any rank-$k$ output after $q$ adaptive products. This inequality applies to arbitrary transcript-dependent outputs, not only to a Krylov subspace. Combined with the positive spectral edge of order $(k/d)^2$, it forces $q=\Omega(d)$ at dimension $d\asymp_p k\eps^{-1/3}$.

Growing $p$ requires more than total overlap. Schatten residuals increasingly emphasize the largest remaining singular values, so we control an initial sequence of small eigenvalues in the output complement. These eigenvalues retain the scale of a fresh square Gaussian factor, whereas the optimal residual sees the shifted rectangular edge. Comparing the resulting sums of spectral powers gives the uniform norm dependence. The fixed-$p$ argument remains useful both as a transparent explanation of the rank factor and as the proof for $1\le p<2$.

\section{Related Work}
\label{sec:related}
\paragraph{Randomized low-rank approximation and Krylov methods.}
Halko, Martinsson, and Tropp develop a modular framework for randomized range finding and matrix decompositions \citep{HMT11}. Musco and Musco show how block Krylov spaces improve the accuracy dependence of spectral approximation \citep{MM15}. BCW sharpen finite Schatten-norm query bounds \citep{BCW22}; Bakshi and Narayanan further remove input-dimension dependence from the fixed-$p$ rank-one upper bound, up to accuracy logarithms \citep{BN23}. Meyer, Musco, and Musco analyze single-vector Krylov methods, including the role of spectral gaps and random perturbations \citep{MMM24}. Chen et al.\ extend block-size comparisons to intermediate blocks \citep{CEMMR26}. Kacham and Woodruff improve arithmetic running times and analyze finite-precision Schatten approximation \citep{KW24}. These results concern algorithms and their implementation; the upper bounds used here are specifically those of BCW and Musco--Musco. We address the information required by any adaptive algorithm, regardless of its update rule.

\paragraph{Adaptive matrix learning and lower bounds.}
Braverman et al.\ use a Gaussian Wishart ensemble to prove oracle lower bounds for linear regression and eigenvalue estimation \citep{BHSW20}. Their adaptive decomposition leaves a fresh Wishart block after conditioning on queries and replies. BCW use a shifted square Wishart matrix to reduce low-rank approximation to small-eigenvalue estimation \citep{BCW22}. Our proof retains these ideas but uses a rectangular factor with exactly $k$ missing columns. The new overlap analysis tracks both posterior covariance outside the query span and the nullspace mass accumulated inside it; this is what yields a rank-dependent information bound. Amsel et al.\ study approximation from more general structured matrix families \citep{AmselEtAl26}. Their framework situates low-rank approximation within matrix learning; its low-rank comparison concerns Frobenius error. Kacham and Woodruff prove measurement--round tradeoffs for adaptive general linear sensing, including constant-factor Schatten approximation \citep{KW23}. That result uses a different measurement budget and does not give the joint rank--accuracy lower bound for sequential exact matvec queries established here.

\paragraph{Random-matrix estimates.}
Rudelson and Vershynin provide rectangular smallest-singular-value estimates \citep{RV09}, and Wei controls intermediate singular values \citep{Wei17}. We use these results to compare rectangular and square spectral profiles after adaptive queries. At bounded rank, coarse profile constants do not separate the relevant eigenvalues. Edelman computes smallest-Wishart-eigenvalue laws \citep{Edelman88,Edelman91}, and Nagao--Forrester give a broader hard-edge distributional analysis \citep{NF98}. We use only the square limiting law and the one-extra-row formula, with their normalization made explicit below, to obtain the probability margin for success probability $2/3$. The term hard edge refers here to the smallest eigenvalues near zero. These distributional inputs are inherited; their use for an arbitrary adaptive output complement is established below.

\section{Oracle Model and Main Results}
\label{sec:model}
The public parameters are $m,n,k,p,\eps$, with $N=\min\{m,n\}$ and $1\le k<N$. An input is any real $m\times n$ matrix. At each step a measurable algorithm chooses either $v\in\R^n$ and receives $Av$, or $u\in\R^m$ and receives $A^\top u$. One vector product costs one query; a block of $b$ vectors costs $b$ queries. Replies are exact, and computation and storage are uncharged. The private random seed is independent of the input. An output is any transcript-measurable rank-$k$ orthogonal projector $P\in\R^{n\times n}$, with success defined by
\begin{equation}
 \|A(I-P)\|_{S_p}\le(1+\eps)\OPT_{p,k}(A),\qquad
 \OPT_{p,k}(A)=\min_{\rank(P_*)=k}\|A(I-P_*)\|_{S_p}.
 \label{eq:model}
\end{equation}
The minimum is over orthogonal projectors. For finite $p$, the Schatten norm is the $\ell_p$ norm of the singular values; at $p=\infty$, it is their maximum.

Define $Q_p^*(m,n,k,\eps)$ as the least integer $q$ for which an algorithm makes at most $q$ queries on every input and every realization of its random seed and satisfies \eqref{eq:model} with probability at least $2/3$ on every input. Thus the budget is a worst-case upper limit, not an expected query count. A policy that stops early can be padded with uninformative calls. All probability thresholds below refer to this fixed convention; we do not assert optimal dependence on a variable confidence parameter.

For each fixed input size, all random variables used in conditional-law arguments take values in standard Borel spaces. We include the seed in the transcript sigma-field. Regular conditional distributions are asserted for almost every realized transcript under its induced law. We write $\Pi_V$ for the orthogonal projector onto a subspace $V$ and $\lambda_j^\uparrow$ for increasingly ordered eigenvalues. Tildes in complexity statements suppress polylogarithmic factors in dimension, rank, norm parameter, and inverse accuracy; in fixed-$p$ statements constants may depend on $p$.
\begin{theorem}[Uniform norm dependence and the spectral endpoint]
\label[theorem]{thm:main}
There exist universal constants $C,\eps_0>0$ such that the following holds.
Let
\[
 2\le p<\infty,
 \qquad 0<\eps\le\eps_0,
 \qquad 1\le k<N,
\]
and define
\begin{equation}
 \Phi(p,\eps):=\eps^{-1/3}\min\{p^{1/6},\eps^{-1/6}\}.
 \label{eq:Phi}
\end{equation}
Then every randomized fully adaptive exact two-sided matvec algorithm for
\eqref{eq:model} requires
\begin{equation}
 {
 Q_p^*(m,n,k,\eps)
 \ge C^{-1}\min\{N,k\Phi(p,\eps)\}.}
 \label{eq:main-lower}
\end{equation}
Consequently,
\begin{equation}
 {
 Q_p^*(m,n,k,\eps)
 =\wtTheta\!\left(\min\{N,k\Phi(p,\eps)\}\right).}
 \label{eq:main-matching}
\end{equation}
For $p=\infty$, the corresponding endpoint is
\begin{equation}
 {
 Q_\infty^*(m,n,k,\eps)
 =\wtTheta\!\left(
     \min\left\{N,\frac{k}{\sqrt\eps}\right\}
   \right).}
 \label{eq:spectral-main}
\end{equation}
We restrict to $k<N$. If $n\le m$ and $k=N=n$, then $P=I_n$ is a zero-query solution. When $m<n$ and $k=N=m$, the output must instead identify a right subspace containing the row space; that endpoint is not claimed here.
\end{theorem}

The main theorem has constants independent of $p$. It therefore applies when $p$ grows with dimension or inverse accuracy. The range below two is supplied by the following complementary result; its constants are allowed to depend on the chosen norm.

\begin{proposition}[Fixed norms below two]
\label[proposition]{prop:small-p}
Fix $1\le p<2$. There exist $c_p,C_p,\eps_p>0$ such that, for every $0<\eps\le\eps_p$ and $1\le k<N$,
\begin{equation}
 c_p\min\{N,k\eps^{-1/3}\}
 \le Q_p^*(m,n,k,\eps)
 \le C_p\min\{N,k\eps^{-1/3}\}\operatorname{polylog}(N,k,1/\eps).
 \label{eq:small-p-rate}
\end{equation}
The oracle, algorithm class, and success criterion are those of \Cref{thm:main}.
\end{proposition}

The lower bound follows from the posterior-overlap argument in \Cref{sec:fixed}; the upper bound is inherited from \citet{BCW22} and exact recovery. This argument also provides a direct proof of the fixed-norm consequence below.

\begin{corollary}[Every fixed finite Schatten norm]
\label[corollary]{cor:fixed}
For each fixed $1\le p<\infty$, every $1\le k<N$, and sufficiently small $\eps$ (with threshold allowed to depend on $p$),
\begin{equation}
 Q_p^*(m,n,k,\eps)=\wtTheta_p\!\left(\min\{N,k\eps^{-1/3}\}\right).
 \label{eq:fixed-rate}
\end{equation}
\end{corollary}
\begin{proof}
For $1\le p<2$, apply \Cref{prop:small-p}. For fixed $p\ge2$, restrict to $\eps\le\min\{\eps_0,1/p\}$. Then $p\eps\le1$, so $\Phi(p,\eps)=p^{1/6}\eps^{-1/3}$ in \Cref{thm:main}. For any fixed $a\ge1$ and $x\ge0$, $\min\{N,x\}\le\min\{N,ax\}\le a\min\{N,x\}$. Absorbing $p^{1/6}$ and fixed-$p$ logarithms into the $p$-dependent constants proves the claim, including the dimension cap.
\end{proof}

Thus the fixed-norm corollary does not assert uniformity along a varying sequence $p=p(\eps)$. The uniform theorem supplies that information for $p\ge2$. Its transition occurs at $p\eps\asymp1$ in polynomial order; the corridor from $1/\eps$ to $(\log N)/\eps$ changes only logarithmic factors.

\begin{corollary}[Near-full-rank saturation]
\label[corollary]{cor:near-full}
Let $k=N-r$, where $1\le r<N$, and use the accuracy range of
\Cref{thm:main}. For $2\le p<\infty$,
\[
 Q_p^*(m,n,N-r,\eps)
 =\wtTheta\!\left(\min\{N,(N-r)\Phi(p,\eps)\}\right).
\]
For $p=\infty$, the corresponding rate is
$\wtTheta(\min\{N,(N-r)/\sqrt\eps\})$.
Thus the finite-$p$ rate is linear up to logarithms when
$(N-r)\Phi(p,\eps)\ge N$, and the spectral rate is linear when
$N-r\ge N\sqrt\eps$. In particular, fixed codimension and fixed admissible
accuracy give linear complexity as $N\to\infty$.
\end{corollary}
\begin{proof}
Substitute $k=N-r$ in \Cref{thm:main}.
\end{proof}

\subsection{Why the rank and accuracy costs multiply}
For the fixed-$p$ proof, set $R=\min\{N,k\eps^{-1/3}\}$; for the uniform theorem, use $R=\min\{N,k\Phi\}$ or its spectral counterpart. There are two branches. If $R\le Mk$, with a sufficiently large constant $M$ (depending on $p$ only in the fixed-$p$ argument), exact recovery of a generic rank-$k$ Gaussian range gives the required lower bound. Otherwise, take $d=\lfloor\eta R\rfloor\gg k$ and use the rectangular Wishart family. This second branch includes $R=N\gg k$. Reaching the ambient cap alone does not justify the rank baseline.

The shared information statement is \Cref{lem:deferred}. Its induction conditions on an explicitly enlarged sigma-field while leaving the remaining Gaussian block independent. The enlargement is a proof device and is removed before making claims about the algorithm. \Cref{thm:posterior-overlap} then bounds the overlap of an arbitrary output with the hidden nullspace. On a high-probability spectral event, the approximation criterion forces an overlap of at least $3k/4$ whenever $d\lesssim_p k\eps^{-1/3}$. This proves \Cref{prop:small-p} and gives a direct proof of \Cref{cor:fixed}.

For \Cref{thm:main}, \Cref{thm:intermediate} shows that the output complement retains small eigenvalues of scale $j^2/d^2$, even after adaptive queries. In contrast, the positive eigenvalues of the full Wishart matrix have scale at least $(k+j)^2/d^2$. The residual power is a sum of terms $(1-\theta_j)^{p/2}$. At $pk^2/d^2=O(1)$, a comparison of $\Theta(k)$ modes gives excess of order $pk^3/d^2$, while the optimal residual power is at most a constant times $d/\sqrt p$. Relative accuracy therefore fails at $d^3\lesssim k^3\sqrt p/\eps$. At larger $pk^2/d^2$, an exponential comparison yields the spectral scale $d\lesssim k/\sqrt\eps$. Bounded $k$ requires an exact one-mode separation; its probability is kept strictly above $1/3$ in \Cref{lem:finite-k-separation}.

\subsection{Embedding the hard dimension}
The lower bounds are proved on symmetric matrices, so transpose queries give no additional information. The following reduction embeds them into the rectangular model while allowing an arbitrary output subspace.

\begin{lemma}[Zero-padding monotonicity]
\label[lemma]{lem:padding}
Let $k\le d\le\min\{m,n\}$. A lower bound for $d\times d$ inputs transfers to $m\times n$ inputs by adding zero rows and columns.
\end{lemma}
\begin{proof}
Write $\bar A=J_m A J_n^\top$, where $J_m,J_n$ are the coordinate embeddings from $\R^d$. For a padded output $\bar P=\bar Z\bar Z^\top$, put $Z_1=J_n^\top\bar Z$. Complete $\range(Z_1)$ to a $k$-subspace and let $P$ be its projector. Then $Z_1Z_1^\top\preceq P$, so
\[
 A(I-P)A^\top\preceq A(I-Z_1Z_1^\top)A^\top.
\]
The matrix on the right, embedded by $J_m$, equals $\bar A(I-\bar P)\bar A^\top$. Eigenvalue monotonicity and the identity $(I-P)^2=I-P$ compare the residual norms, including $p=\infty$. Padding leaves the optimal tail singular values unchanged. Each padded query can be simulated with one original product by restricting its input and embedding its reply.
\end{proof}

\section{Adaptive Gaussian Information}
\label{sec:posterior}
An orthogonalized query contains all new information in a symmetric exact-product oracle: its component in the previous query span has a known reply. We therefore count only informative, orthonormal queries. If the algorithm uses fewer than a prescribed number, append orthogonal queries and ignore their replies in its output. The results below apply to this stronger transcript as well.

\subsection{A rectangular deferred-decision theorem}
The square-Wishart conditioning argument of \citet[Lemma~13 and Appendix~C.2]{BHSW20} is the starting point. We need a rectangular version because $X\in\R^{d\times(d-k)}$ leaves an exact $k$-dimensional nullspace. The strengthened induction below makes both the latent rotations and the conditioning explicit.

\begin{lemma}[Rectangular adaptive Gaussian factorization]
\label[lemma]{lem:deferred}
Let $1\le s<d$, let $X\in\R^{d\times s}$ have independent $N(0,1/d)$ entries, and put $W=XX^\top$. Consider a measurable adaptive policy with independent seed $\xi$, orthonormal queries $v_1,\ldots,v_q$, and $q\le s$. Set
\[
 \cF_t=\sigma(\xi,v_1,Wv_1,\ldots,v_t,Wv_t),\qquad \cF_0=\sigma(\xi).
\]
For every $0\le t\le q$, there are an $\cF_t$-measurable orthogonal row matrix $U_t$ whose first $t$ rows are $v_1^\top,\ldots,v_t^\top$, an enlarged sigma-field $\cG_t\supseteq\cF_t$, and a $\cG_t$-measurable orthogonal matrix $O_t\in\R^{s\times s}$ such that
\begin{equation}
 U_tXO_t=\begin{bmatrix}A_t&0\\ C_t&D_t\end{bmatrix},\quad
 A_t\in\R^{t\times t},\quad
 D_t\in\R^{(d-t)\times(s-t)}.
 \label{eq:block-factor}
\end{equation}
The blocks $A_t,C_t$ are $\cF_t$-measurable, and $A_t$ is lower triangular with positive diagonal almost surely. Moreover, for every bounded Borel function $\varphi$,
\begin{equation}
 \E[\varphi(D_t)\mid\cG_t]=\int\varphi\,d\gamma_t
 \quad\text{almost surely},
 \label{eq:conditional-G}
\end{equation}
where $\gamma_t$ is the fixed product $N(0,1/d)$ law in the displayed dimensions. In particular, $D_t$ is independent of $\cG_t$, and
\begin{equation}
 \E[\varphi(D_t)\mid\cF_t]=\int\varphi\,d\gamma_t
 \quad\text{almost surely}.
 \label{eq:conditional-F}
\end{equation}
Thus its regular conditional law is Gaussian for almost every realized query transcript. No assertion is made about conditional-law versions on null sets of transcripts.
\end{lemma}

\begin{proof}
We prove all statements simultaneously, including \eqref{eq:conditional-G}, and construct $\cG_t$ recursively. Initially let $U_0=I_d$, $O_0=I_s$, $\cG_0=\sigma(\xi)$, and $D_0=X$. Independence of $X$ and $\xi$ proves the base case.

Suppose the claims hold at time $t<s$. Since $v_{t+1}$ is $\cF_t$-measurable and orthogonal to the first $t$ queries, choose an $\cF_t$-measurable orthogonal transformation $R_t$ on the last $d-t$ row coordinates that sends its coordinate vector to the first basis vector. Measurable completions can be fixed, for example, by Gram--Schmidt applied to the standard basis with the first nonzero pivot. Set
\[
 U_{t+1}=\operatorname{diag}(I_t,R_t)U_t.
\]
Before changing the right coordinates, write
\[
 U_{t+1}XO_t=
 \begin{bmatrix}
 A_t&0\\ c_t^\top&z_t^\top\\ C_t'&D_t'
 \end{bmatrix}.
\]
Here $c_t,C_t'$ are $\cF_t$-measurable. Conditional on $\cG_t$, the row $z_t^\top$ and the entries of $D_t'$ are mutually independent $N(0,1/d)$ variables. This follows from the induction hypothesis and left orthogonal invariance, because $R_t$ is already $\cG_t$-measurable.

Put $\rho_t=\|z_t\|_2$ and $u_t=z_t/\rho_t$; $\rho_t>0$ almost surely. Choose a Borel orthogonal completion $[u_t,V_t]$ and set $b_t=D_t'u_t$. The transformed reply to $v_{t+1}$ is
\begin{equation}
 U_{t+1}Wv_{t+1}=
 \begin{bmatrix}
 A_tc_t\\ \|c_t\|_2^2+\rho_t^2\\ C_t'c_t+\rho_tb_t
 \end{bmatrix}.
 \label{eq:reply-decomposition}
\end{equation}
It follows that, modulo null sets,
\[
 \cF_{t+1}=\cF_t\vee\sigma(\rho_t,b_t).
\]
Define the enlarged field and right rotation by
\begin{equation}
 \cG_{t+1}:=\cG_t\vee\sigma(z_t,b_t)
 =\cG_t\vee\cF_{t+1}\vee\sigma(u_t),\qquad
 O_{t+1}:=O_t\operatorname{diag}(I_t,[u_t,V_t]).
 \label{eq:enlarged-filtration}
\end{equation}
In particular, $\cG_{t+1}$ contains the old enlarged field and the new actual transcript, but not the unobserved orthogonal Gaussian block. Direct multiplication gives
\[
 A_{t+1}=\begin{bmatrix}A_t&0\\c_t^\top&\rho_t\end{bmatrix},
 \qquad C_{t+1}=[C_t'\ b_t],\qquad
 D_{t+1}=D_t'V_t.
\]
The first two blocks are $\cF_{t+1}$-measurable, and $A_{t+1}$ is lower triangular with positive diagonal.

Conditional on $(\cG_t,z_t)$, orthogonal Gaussian decomposition of each row of $D_t'$ shows that $D_t'u_t$ and $D_t'V_t$ are independent. The latter has the same product Gaussian law for every value of $(\cG_t,z_t)$. Conditioning further on $b_t=D_t'u_t$ therefore leaves this law unchanged. This proves
\[
 \E[\varphi(D_{t+1})\mid\cG_t\vee\sigma(z_t,b_t)]
 =\int\varphi\,d\gamma_{t+1},
\]
which is the strengthened induction claim. The construction also makes $O_{t+1}$ measurable with respect to $\cG_{t+1}$, completing the induction. At $t=s$, the hidden block has zero columns and its law is the corresponding point mass.

Finally, the tower property yields \eqref{eq:conditional-F}. Equivalently, a regular conditional kernel for $D_t$ given the transcript equals $\gamma_t$ outside a set of transcript probability zero. The exposed leading block is automatically the positive-diagonal Cholesky factor of the queried Gram matrix; no further latent conditioning is needed.
\end{proof}

\subsection{Exact recovery as the rank baseline}
The same lemma shows why exactly rank-$k$ inputs require $k$ products. This observation will be used only when the target lower bound is at most a constant times $k$.
\begin{lemma}[All-rank exact range-recovery baseline]
\label[lemma]{lem:range-baseline}
Let $1\le k<d$, draw
\[
 G\in\R^{d\times k},
 \qquad G_{ij}\overset{\mathrm{iid}}\sim N(0,1/d),
 \qquad A:=GG^\top,
\]
and give an algorithm exact adaptive matvec access to $A$.  Any randomized
algorithm that outputs a rank-$k$ projector $P$ satisfying
\begin{equation}
 A(I-P)=0
 \label{eq:exact-range-goal}
\end{equation}
with any fixed positive success probability requires at least $k$ queries.
Consequently, for every $1\le p\le\infty$ and every $\eps>0$, rank-$k$
relative-error Schatten-$p$ LRA has a $k$-query lower bound for every
nontrivial rank $1\le k<d$.
\end{lemma}

\begin{proof}
By Yao's principle, fix a deterministic algorithm making $q<k$ adaptive
queries.  The case $q=0$ is immediate from the nonatomic distribution of
$\range(G)$, so assume $q\ge1$.  Apply \Cref{lem:deferred} with
column count $s=k$.  Conditional on the actual transcript, in a
transcript-measurable row basis
and an invisible right-orthogonal gauge,
\[
 G_q=
 \begin{bmatrix}
 B_q&0\\
 C_q&D_q
 \end{bmatrix},
 \qquad
 D_q\in\R^{(d-q)\times(k-q)},
\]
where $B_q$ is invertible almost surely and $D_q$ is an independent iid
Gaussian matrix.  Let
\[
 E_q:=\range\!\begin{bmatrix}B_q\\C_q\end{bmatrix}
\]
be the exposed $q$-dimensional part of $\range(G)$.  Define the
transcript-measurable quotient map
\[
 T_q:\R^q\oplus\R^{d-q}\to\R^{d-q},
 \qquad
 T_q(x,y):=y-C_qB_q^{-1}x.
\]
Then $\ker(T_q)=E_q$ and
\begin{equation}
 T_q(\range(G))=\range(D_q).
 \label{eq:quotient-hidden-range}
\end{equation}

The output subspace $\range(P)$ is fixed conditional on the transcript.  If
\eqref{eq:exact-range-goal} holds, then
$\range(I-P)\subseteq\ker(A)$.  Both spaces have dimension $d-k$, so they are
equal, and therefore
\[
 \range(P)=\range(A)=\range(G).
\]
Equation~\eqref{eq:quotient-hidden-range} would then force the fixed subspace
$T_q(\range(P))$ to equal $\range(D_q)$.  But
\[
 1\le k-q<d-q
\]
and the column space of an iid Gaussian
$(d-q)\times(k-q)$ matrix has a nonatomic Haar distribution on
$\mathrm{Gr}(k-q,d-q)$.  Hence this equality has conditional probability
zero.  Every deterministic $q<k$ algorithm therefore succeeds with
probability zero under the Gaussian input distribution, and Yao's principle
gives the claimed randomized lower bound.

Finally, $A$ has rank exactly $k$, so the optimal rank-$k$ residual is zero in
every Schatten norm.  Relative-error success therefore implies
\eqref{eq:exact-range-goal}.  Since $A$ is symmetric positive semidefinite,
the result applies to the full two-sided oracle.
\end{proof}

\begin{remark}[Sharpness on exactly rank-$k$ inputs]
The baseline is exact: $k$ independent Gaussian query vectors form a matrix
$\Omega\in\R^{d\times k}$ for which $G^\top\Omega$ is invertible almost
surely, so the responses
$A\Omega=G(G^\top\Omega)$ span $\range(G)$.  Thus generic rank-$k$ PSD inputs
have exact range-recovery complexity $k$.  The full LRA upper bound remains
capped by $N$ because arbitrary inputs need not have rank $k$.
\end{remark}

\subsection{Posterior nullspace overlap}

Specialize \Cref{lem:deferred} to
$X\in\R^{d\times(d-k)}$ and $W=XX^\top$.  Put
$K:=\ker W$.  After $q$ queries, write $s=d-q$ and use the block
factorization~\eqref{eq:block-factor}.  A vector $(x,y)\in\R^q\oplus\R^s$
lies in $K$ exactly when
\[
  A_q^\top x+C_q^\top y=0,
  \qquad
  D_q^\top y=0.
\]
Thus, with
\[
  B_q:=A_q^{-\top}C_q^\top,
  \qquad
  K_q^0:=\ker(D_q^\top),
\]
we have the graph representation
\begin{equation}
  K=\left\{\binom{-B_qy}{y}:y\in K_q^0\right\}.
  \label{eq:graph-K}
\end{equation}
For almost every transcript, $K_q^0$ has the Haar law of a $k$-subspace of
$\R^s$.

\begin{lemma}[Bottom-block posterior covariance]
\label[lemma]{lem:bottom-covariance}
Let $\overline\Pi_q:=\E[\Pi_K\mid\cF_q]$ in the query-adapted basis and write
\[
  \overline\Pi_q=
  \begin{bmatrix}
    M_{11}&M_{12}\\
    M_{12}^\top&M_{22}
  \end{bmatrix}
  \quad\text{on }\R^q\oplus\R^s.
\]
Then
\begin{equation}
  0\preceq M_{22}\preceq\frac{k}{s}I_s.
  \label{eq:bottom-bound}
\end{equation}
\end{lemma}

\begin{proof}
Let $Q_0\in\R^{s\times k}$ have orthonormal columns spanning $K_q^0$ and set
$T=[-B_q^\top,I_s]^\top$.  The graph projector is
\[
  \Pi_K=TQ_0(Q_0^\top T^\top TQ_0)^{-1}Q_0^\top T^\top.
\]
Its bottom-right block is
\[
  Q_0\bigl(Q_0^\top(I_s+B_q^\top B_q)Q_0\bigr)^{-1}Q_0^\top
  \preceq Q_0Q_0^\top.
\]
Taking conditional expectation and using
$\E[Q_0Q_0^\top]=(k/s)I_s$ proves~\eqref{eq:bottom-bound}.
\end{proof}

\begin{lemma}[Nullspace mass learned by the query span]
\label[lemma]{lem:query-mass}
For every adaptive orthonormal query sequence with span $V_q$,
\begin{equation}
  \E\tr(\Pi_{V_q}\Pi_K)
  \le k\sum_{t=0}^{q-1}\frac{1}{d-t}
  \le\frac{kq}{d-q}.
  \label{eq:query-mass}
\end{equation}
\end{lemma}

\begin{proof}
Write the orthonormal queries as $v_1,\ldots,v_q$.  Conditional on
$\cF_{t-1}$, the vector $v_t$ is a unit vector in the hidden row block of
dimension $d-t+1$.  \Cref{lem:bottom-covariance} at time $t-1$ gives
\[
  \E[v_t^\top\Pi_Kv_t\mid\cF_{t-1}]
  \le\frac{k}{d-t+1}.
\]
Summing and applying the tower property proves the claim.
\end{proof}

\begin{theorem}[Adaptive posterior overlap]
\label[theorem]{thm:posterior-overlap}
After any $q\le d-k$ randomized, fully adaptive products with $W$, every
transcript-measurable rank-$k$ projector $P$ satisfies
\begin{equation}
  {
  \E\tr(P\Pi_K)
  \le\frac{2k(q+k)}{d-q}.}
  \label{eq:posterior-overlap}
\end{equation}
Consequently, if $q,k\le\delta d$ for a sufficiently small universal
constant $\delta>0$, then
\begin{equation}
  \Prb\!\left[\tr(P\Pi_K)\ge\frac{3k}{4}\right]\le0.1.
  \label{eq:overlap-probability}
\end{equation}
\end{theorem}

\begin{proof}
Condition on the final transcript.  Every positive-semidefinite block matrix
satisfies
\[
  \begin{bmatrix}A&C\\C^\top&D\end{bmatrix}
  \preceq
  2\begin{bmatrix}A&0\\0&D\end{bmatrix}.
\]
Applying this to $\overline\Pi_q$ and using
\Cref{lem:bottom-covariance},
\begin{align*}
  \E[\tr(P\Pi_K)\mid\cF_q]
  &=\tr(P\overline\Pi_q)\\
  &\le2\tr(M_{11})+2k\|M_{22}\|_{\op}\\
  &\le2\tr(\Pi_{V_q}\overline\Pi_q)+\frac{2k^2}{d-q}.
\end{align*}
Take expectations and apply \Cref{lem:query-mass}.  Markov's inequality
gives
\[
  \Prb\!\left[\tr(P\Pi_K)\ge\frac{3k}{4}\right]
  \le\frac{8}{3}\frac{q+k}{d-q},
\]
which is at most $0.1$ for a sufficiently small $\delta$.
\end{proof}

\section{The Fixed-Norm Lower Bound}\label{sec:fixed}
We first turn missing nullspace information into Schatten error with constants allowed to depend on $p$. Conditioning the spectrum away from zero makes the argument valid even when $p/2<1$.

\begin{lemma}[Well-conditioned Schatten compression]
\label[lemma]{lem:compression}
Let $H\in\R^{d\times d}$ be symmetric positive definite with eigenvalues in a
fixed interval $[a,b]\subset(0,\infty)$.  Let $P$ be any rank-$k$ projector
and define
\[
  \Delta_H(P)
  :=\sum_{i=1}^k\lambda_i(H)^2-\tr(PH^2),
\]
where the eigenvalues of $H$ are ordered decreasingly.  For fixed
$1\le p<\infty$, there is $c_{p,a,b}>0$ such that
\begin{equation}
  \|H(I-P)\|_{S_p}^p
  -\min_{\rank(P_*)=k}\|H(I-P_*)\|_{S_p}^p
  \ge c_{p,a,b}\Delta_H(P).
  \label{eq:compression}
\end{equation}
Consequently, if $P$ satisfies~\eqref{eq:model} and $\eps\le1$, then
\begin{equation}
  \Delta_H(P)\le C_{p,a,b}\eps d.
  \label{eq:compression-consequence}
\end{equation}
\end{lemma}

\begin{proof}
Set $C=H^2$ and write its eigenvalues as
$c_1\ge\cdots\ge c_d\in[a^2,b^2]$.  Put $Q=I-P$, and let
$\theta_1\ge\cdots\ge\theta_{d-k}$ be the eigenvalues of $QCQ$ on
$\range(Q)$.  Poincare separation gives $\theta_j\ge c_{k+j}$.  With
$r=p/2$, the derivative of $x^r$ has a positive lower bound on
$[a^2,b^2]$, so
\[
  \theta_j^r-c_{k+j}^r
  \ge c_{p,a,b}(\theta_j-c_{k+j}).
\]
Summing and using
\[
  \sum_{j=1}^{d-k}\theta_j-\sum_{j=k+1}^dc_j
  =\sum_{i=1}^kc_i-\tr(PC)
  =\Delta_H(P)
\]
proves~\eqref{eq:compression}.  Under~\eqref{eq:model}, the excess in the
left-hand side is at most
$((1+\eps)^p-1)\OPT^p\le C_{p,b}\eps d$.
\end{proof}

\subsection{A well-conditioned hard family}

Draw
\begin{equation}
  X\in\R^{d\times(d-k)},
  \qquad
  X_{ij}\overset{\mathrm{iid}}\sim N(0,1/d),
  \qquad
  W:=XX^\top.
  \label{eq:wishart}
\end{equation}
Then $K=\ker W$ is Haar distributed on the Grassmannian of $k$-subspaces.

\begin{lemma}[Wishart spectral event]
\label[lemma]{lem:spectral-event}
There are universal constants $c_0,C_0,C>0$ such that, for $d\ge Ck$,
\begin{equation}
  \Prb\!\left[
    \|W\|_{\op}\le C_0
    \text{ and }
    \lambda_{\min}^+(W)\ge c_0(k/d)^2
  \right]
  \ge0.95.
  \label{eq:spectral-event}
\end{equation}
\end{lemma}

\begin{proof}
The operator-norm bound is standard for Gaussian matrices.  For the lower
edge, the smallest-singular-value estimate of Rudelson and
Vershynin~\citep[Theorem~1.1]{RV09} gives
$s_{\min}(X)\ge c k/d$ with arbitrarily high fixed probability because the
unscaled rectangular edge is
$\sqrt d-\sqrt{d-k-1}\asymp k/\sqrt d$.  Squaring proves the claim.
\end{proof}

Fix $L:=2C_0$ and set
\begin{equation}
  H:=I-W/L.
  \label{eq:H-hard}
\end{equation}
On the event in \Cref{lem:spectral-event}, $H$ has spectrum in
$[1/2,1]$ and top-$k$ eigenspace $K$.

\begin{lemma}[Relative Schatten error forces nullspace overlap]
\label[lemma]{lem:lra-to-overlap}
On the event~\eqref{eq:spectral-event}, if a rank-$k$ projector $P$ satisfies
\eqref{eq:model} for $H$, then
\begin{equation}
  k-\tr(P\Pi_K)
  \le C_p\frac{\eps d^3}{k^2}.
  \label{eq:overlap-deficit}
\end{equation}
In particular, if $d\le\eta_p k\eps^{-1/3}$ for sufficiently small
$\eta_p>0$, then $\tr(P\Pi_K)\ge3k/4$.
\end{lemma}

\begin{proof}
The top $k$ eigenvalues of $H$ are one.  \Cref{lem:compression} gives
$k-\tr(PH^2)\le C_p\eps d$.  Since
\[
  I-H^2=\frac{2}{L}W-\frac{1}{L^2}W^2
  \succeq\frac{3}{2L}W
\]
on $\|W\|_{\op}\le L/2$, we obtain $\tr(PW)\le C_p'\eps d$.  Also
\[
  W\succeq c_0(k/d)^2(I-\Pi_K),
\]
so
\[
  \tr(PW)
  \ge c_0(k/d)^2\bigl(k-\tr(P\Pi_K)\bigr).
\]
Combining the two inequalities proves~\eqref{eq:overlap-deficit}.
\end{proof}

A product with $H$ is equivalent to one product with $W$ because
$Wv=L(v-Hv)$.  Since $H$ is symmetric, two-sided access gives no additional
power on this hard family.

\subsection{Completing the fixed-norm bound}

\begin{proof}[Posterior-overlap proof for fixed finite $p$]
Fix $1\le p<\infty$. We prove the lower bound directly throughout this fixed-norm range; its restriction to $p<2$ supplies \Cref{prop:small-p}.
Set
\[
  R:=\min\{N,k\eps^{-1/3}\}.
\]
If $R\le M_pk$ for a sufficiently large constant $M_p$, use the Gaussian
rank-$k$ hard instance from \Cref{lem:range-baseline} in dimension $N$
(or in an $N\times N$ block of the prescribed rectangular dimensions).  It
gives
\[
  Q_p^*(m,n,k,\eps)\ge k\ge R/M_p.
\]

Assume $R>M_pk$.  Choose a sufficiently small constant $\eta_p>0$ and let
$d=\lfloor\eta_pR\rfloor$.  Taking $M_p$ large after $\eta_p$ is fixed
ensures
\[
  d\le N,
  \qquad
  d\le\eta_pk\eps^{-1/3},
  \qquad
  k\le d/100.
\]
Construct the $d\times d$ hard matrix~\eqref{eq:H-hard} and zero-pad it using
\Cref{lem:padding}.
Increase $M_p$ further so that $d\ge Ck$ for the constant in
\Cref{lem:spectral-event} and $d\ge\eta_pR/2$.

Suppose an algorithm uses $q\le d/100$ queries and outputs $P$.  On the
spectral event of \Cref{lem:spectral-event}, success implies
$\tr(P\Pi_K)\ge3k/4$ by \Cref{lem:lra-to-overlap}.  But
\Cref{thm:posterior-overlap} gives
\[
  \Prb\!\left[\tr(P\Pi_K)\ge3k/4\right]\le0.1.
\]
Without conditioning the hard distribution on the spectral event, the union bound gives success probability at most $0.1+0.05=0.15$. The hard distribution is independent of the algorithm's private seed, and the posterior bound already includes that seed. A $2/3$-successful algorithm on every input would also have distributional success at least $2/3$, a contradiction. Hence its worst-case budget is $\Omega_p(d)=\Omega_p(R)$.
\end{proof}

Combining this lower bound with the upper bound of~\citep{BCW22} and exact
recovery proves \Cref{prop:small-p} and gives an independent proof of \Cref{cor:fixed}.

\section{Spectral Information Beyond Total Overlap}\label{sec:spectrum}
\subsection{Rectangular spectral profile}

We now use a universal normalization, independent of $p$, for the growing-norm analysis. Fix a hard dimension $d$ and draw
\begin{equation}
 X\in\R^{d\times(d-k)},
 \qquad X_{ij}\overset{\mathrm{iid}}\sim N(0,1/d),
 \qquad W:=XX^\top.
 \label{eq:XW}
\end{equation}
Let $L$ be a sufficiently large universal constant and set
\begin{equation}
 H:=I-W/L,
 \qquad S:=I-H^2=\frac{2}{L}W-\frac{1}{L^2}W^2.
 \label{eq:HS}
\end{equation}
On the event $\|W\|_{\op}\le L/2$,
\begin{equation}
 0\preceq S\preceq \frac{2}{L}W,
 \qquad
 S\succeq \frac{3}{2L}W,
 \label{eq:S-W-comparison}
\end{equation}
and $H$ has $k$ eigenvalues equal to one, with all remaining eigenvalues in
$[1/2,1)$.

We order the positive eigenvalues of $W$ increasingly:
\[
 0<\lambda_1^+(W)\le\cdots\le\lambda_{d-k}^+(W),
\]
and likewise write
$0<s_1\le\cdots\le s_{d-k}$ for the positive eigenvalues of $S$.

\begin{lemma}[Uniform rectangular lower spectrum]
\label[lemma]{lem:rect-spectrum}
There are universal constants $c_0,L>0$ and $d_0\in\mathbb N$ such that for
every $d\ge d_0$ and every $1\le k\le d/4$, with probability at least $0.995$,
\begin{equation}
 \|W\|_{\op}\le L/2,
 \qquad
 \lambda_j^+(W)\ge c_0\left(\frac{k+j}{d}\right)^2
 \quad(1\le j\le d-k).
 \label{eq:rect-spectrum}
\end{equation}
Consequently,
\begin{equation}
 s_j\ge c_1\left(\frac{k+j}{d}\right)^2
 \qquad(1\le j\le d-k)
 \label{eq:S-spectrum}
\end{equation}
for a universal $c_1>0$.  In particular, this event is uniform down to
$k=1$.
\end{lemma}

\begin{proof}
Let $r=d-k$.  For each $1\le j\le r$, let $X_j$ be the submatrix formed by
any fixed $r-j+1$ columns of $X$.  Cauchy interlacing for the principal
submatrix $X_j^\top X_j$ gives
\[
 s_j^{\up}(X)\ge s_{\min}(X_j).
\]
Write $G_j=\sqrt d\,X_j$.  It is a $d\times(r-j+1)$ standard Gaussian matrix,
and its row excess plus one is $k+j$.  The rectangular smallest-singular-value estimate of
Rudelson--Vershynin \citep[Theorem~1.1]{RV09} implies that for every
$t\in(0,t_0)$,
\[
 \Prb\!\left[
 s_{\min}(G_j)
 \le t\bigl(\sqrt d-\sqrt{r-j}\bigr)
 \right]
 \le (C t)^{k+j}+e^{-cd}.
\]
Since
\[
 \sqrt d-\sqrt{r-j}
 =\frac{k+j}{\sqrt d+\sqrt{r-j}}
 \ge\frac{k+j}{2\sqrt d},
\]
we obtain
\[
 \Prb\!\left[
 s_j^{\up}(X)\le \frac{t(k+j)}{2d}
 \right]
 \le(Ct)^{k+j}+e^{-cd}.
\]
Choose $t$ so small that
$\sum_{\ell\ge2}(Ct)^\ell<10^{-3}$, and then choose $d_0$ so that
$d e^{-cd}<10^{-3}$.  A union bound over $j$ proves the lower-spectrum
statement with probability at least $0.998$.  For completeness, $1/4$-nets of the two unit spheres have at most
$9^d$ and $9^r$ points. Since each fixed bilinear form of $X$ is
$N(0,1/d)$, the net approximation and a union bound give
\[
 \Prb[\|X\|_{\op}>2t]\le 2\,9^{d+r}e^{-dt^2/2}.
\]
A sufficiently large universal $t$, and then $L\ge8t^2$, make the
operator-norm failure probability at most $0.003$.  Finally,
\[
 S=\frac{2}{L}W-\frac{1}{L^2}W^2
 \succeq\frac{3}{2L}W
\]
on $\|W\|_{\op}\le L/2$, which proves \eqref{eq:S-spectrum}.
\end{proof}

\begin{lemma}[Spectral residual certificate]
\label[lemma]{lem:spectral-certificate}
Let $\mu=\lambda_{\min}^+(W)$ and let $P$ be a rank-$k$ projector, with
$Q=I-P$.  Put
\[
 \delta:=\mu/L,
 \qquad
 \nu(P):=\lambda_{\min}\!\left(
     QWQ\big|_{\range(Q)}
   \right).
\]
Assume $0\preceq W\preceq(L/2)I$, $\eps\le1$, and
\[
 \|HQ\|_{\op}\le(1+\eps)(1-\delta).
\]
Then
\begin{equation}
 \nu(P)
 \ge
 \mu\left(
       1-\frac{\delta}{2}-\frac{3\eps}{2\delta}
     \right).
 \label{eq:spectral-certificate-general}
\end{equation}
In particular:
\begin{enumerate}[label=(\roman*),leftmargin=2em]
\item if $\eps\le\delta/8$, then
\begin{equation}
 \nu(P)\ge\frac{9}{16}\mu;
 \label{eq:spectral-certificate-coarse}
\end{equation}
\item if $\delta\le\eta$ and $\eps\le\eta\delta$ for some
$0<\eta\le1/10$, then
\begin{equation}
 \nu(P)\ge(1-2\eta)\mu.
 \label{eq:spectral-certificate}
\end{equation}
\end{enumerate}
\end{lemma}

\begin{proof}
On $\range(Q)$,
\[
 QH^2Q=I_{\range(Q)}-QSQ,
 \qquad
 \|HQ\|_{\op}^2
 =1-\lambda_{\min}(QSQ|_{\range(Q)}).
\]
Hence success implies
\[
 \lambda_{\min}(QSQ|_{\range(Q)})
 \ge1-(1+\eps)^2(1-\delta)^2.
\]
Since $(1+\eps)^2\le1+3\eps$,
\[
 1-(1+\eps)^2(1-\delta)^2
 \ge2\delta-\delta^2-3\eps.
\]
The comparison $QSQ\preceq(2/L)QWQ$ gives
\[
 \nu(P)
 \ge\frac L2(2\delta-\delta^2-3\eps)
 =\mu\left(
       1-\frac{\delta}{2}-\frac{3\eps}{2\delta}
     \right),
\]
proving \eqref{eq:spectral-certificate-general}.  Because
$\delta\le1/2$ on the assumed operator-norm event, the condition
$\eps\le\delta/8$ yields
$1-\delta/2-3\eps/(2\delta)\ge1-1/4-3/16=9/16$.
The sharper conclusion follows from
$\delta/2+3\eps/(2\delta)\le2\eta$.
\end{proof}

\subsection{Adaptive persistence of small eigenvalues}

The spectral proof needs only the least singular value of the fresh Gaussian
block.  Growing $p$ requires a simultaneous bound for a whole initial segment
of the small singular values.

\begin{theorem}[Adaptive intermediate hard-edge persistence]
\label[theorem]{thm:intermediate}
After $q$ adaptive queries, let $P$ be any $\cF_q$-measurable rank-$k$
projector, set $Q=I-P$, and put
\[
 h:=d-k-q.
\]
If $h\ge1$, then for almost every transcript, conditional on $\cF_q$, with probability at least $0.95$,
all $1\le j\le h$ satisfy
\begin{equation}
 \lambda_j^{\up}\!\left(QWQ\big|_{\range(Q)}\right)
 \le C\frac{j^2}{dh}.
 \label{eq:intermediate-W}
\end{equation}
If $k+q\le d/2$, this becomes
\begin{equation}
 \lambda_j^{\up}\!\left(QWQ\big|_{\range(Q)}\right)
 \le C'\frac{j^2}{d^2}.
 \label{eq:intermediate-W-simple}
\end{equation}
Consequently, if
$\theta_1\le\cdots\le\theta_{d-k}$ are the eigenvalues of
$QSQ|_{\range(Q)}$, then on the same event
\begin{equation}
 \theta_j\le C''\frac{j^2}{d^2}
 \qquad(1\le j\le h)
 \label{eq:intermediate-S}
\end{equation}
whenever $k+q\le d/2$. These eigenvalues are nonnegative when
$\|W\|_{\op}\le L/2$. The conditional probability claim is made before
intersecting with that spectral event.
\end{theorem}

\begin{proof}
Work in the block coordinates of \Cref{lem:deferred}: write
$X_q:=U_qXO_q$ and $W_q:=U_qWU_q^\top=X_qX_q^\top$.
All output-subspace coordinates below use this same row basis. Let
$r=d-k$ and choose $U\in\R^{d\times r}$ with orthonormal columns spanning
$\range(Q)$.  Partition $U=[U_1^\top,U_2^\top]^\top$.  Then
\[
 M:=U^\top X_q=[\,F\ \ U_2^\top D_q\,],
 \qquad F:=U_1^\top A_q+U_2^\top C_q\in\R^{r\times q},
\]
and $U^\top W_qU=MM^\top$. Select a transcript-measurable
$R\in\R^{r\times h}$ with orthonormal columns in $\ker(F^\top)$ and put
$T=U_2R$.  Then $\|T\|_{\op}\le1$ and
\[
 M^\top R=
 \begin{bmatrix}0\\D_q^\top T\end{bmatrix}.
\]
Singular-value interlacing under restriction of the domain gives
\begin{equation}
 s_j^{\up}(M)
 \le s_j^{\up}(D_q^\top T)
 \qquad(1\le j\le h).
 \label{eq:restriction-interlace}
\end{equation}
Conditional on the transcript,
\[
 D_q^\top T\overset{d}{=}d^{-1/2}G_h\Sigma,
 \qquad \Sigma=(T^\top T)^{1/2},
 \qquad \|\Sigma\|_{\op}\le1,
\]
where $G_h$ is standard square Gaussian.  For every ordered singular value,
\[
 s_j^{\up}(G_h\Sigma)
 \le\|\Sigma\|_{\op}s_j^{\up}(G_h)
 \le s_j^{\up}(G_h).
\]
Standard Gaussian entries satisfy the unit-variance, subgaussian, and bounded-density assumptions of \citet[Assumption~1.3]{Wei17}, with fixed universal parameters. Thus \citet[Theorem~1.6]{Wei17} gives
\[
 \Prb[s_j^{\up}(G_h)>C_1tj/\sqrt h]\le e^{-C_2tj}
 \qquad(t>1).
\]
Choose a universal $t$ large enough that $\sum_{j\ge1}e^{-C_2tj}\le0.05$ and apply a union bound. Thus $s_j^{\up}(G_h)\le Cj/\sqrt h$ holds simultaneously for every $j$ with probability at least $0.95$.  Combining this with
\eqref{eq:restriction-interlace} and squaring proves
\eqref{eq:intermediate-W}.  The remaining statements follow from
$QSQ\preceq(2/L)QWQ$.
\end{proof}

The $j=1$ bound supplies the spectral obstruction for sufficiently large $k$. For bounded $k$, its unspecified constants are not enough to compare the fresh and optimal edges. We next retain their exact distributional separation before assembling the spectral result.
\subsection{Bounded-rank edge separation}

The intermediate theorem is sufficient for growing rank, but its coarse
constants do not separate even one mode when $k$ is a small fixed integer.
The exact real-Wishart hard-edge laws provide the strict probability margin needed in that regime.

Let $G_{n+\nu,n}$ be a standard real Gaussian matrix and define
\[
 Y_{n,\nu}:=\sqrt{n+\nu}\,
 s_{\min}(G_{n+\nu,n}).
\]
We need only the square and one-extra-row laws:
\begin{align}
 F_0(t)
 &:=
 \lim_{n\to\infty}\Prb[Y_{n,0}\le t]
 =1-e^{-t-t^2/2},
 \label{eq:F0}\\
 F_1(t)
 &:=
 \lim_{n\to\infty}\Prb[Y_{n,1}\le t]
 =1-e^{-t^2/2}.
 \label{eq:F1}
\end{align}
The square law \eqref{eq:F0} is Edelman's real-Gaussian limit
\citep{Edelman88}; it is also stated, for the scaled least eigenvalue,
in \citet[Lemma~11]{BHSW20}. The change of variable is
$Y_{n,0}^2=n\lambda_{\min}(G_{n,n}^{\top}G_{n,n})$.
For \eqref{eq:F1}, the real Wishart joint eigenvalue density for
$G_{n+1,n}^{\top}G_{n+1,n}$ is proportional to
$e^{-\sum_i\lambda_i/2}\prod_{i<j}|\lambda_i-\lambda_j|$
\citep{Edelman91}. Translating all $n$ eigenvalues by $x$ therefore gives
\[
 \Prb[\lambda_{\min}>x]=e^{-nx/2},\qquad
 \Prb[Y_{n,1}\le t]=1-e^{-nt^2/(2(n+1))}.
\]
Taking $n\to\infty$ gives exactly \eqref{eq:F1}, with the same normalization
as \eqref{eq:F0}. Their strict separation is quantitative. At
\[
 b_0:=\frac34,
 \qquad
 a_0:=\frac45,
\]
we have
\begin{equation}
 F_0(b_0)-F_1(a_0)
 =e^{-8/25}-e^{-33/32}
 >0.36>\frac13.
 \label{eq:CDF-gap}
\end{equation}

\begin{lemma}[Adaptive finite-codimension edge separation]
\label[lemma]{lem:finite-k-separation}
Fix any integer $K\ge2$.  There exist constants
$\alpha_K,\rho_K>0$, $M_K,d_K<\infty$, and
$0<B_K<A_K<A_K^+<\infty$ such that the following holds.
Let
\[
 1\le k<K,
 \qquad
 d\ge d_K,
 \qquad
 d\ge M_K k,
 \qquad
 q\le\alpha_K d.
\]
For every deterministic adaptive $q$-query algorithm and every
$\cF_q$-measurable rank-$k$ output projector $P$, with $Q=I-P$,
\begin{equation}
 \Prb\!\left[
 \begin{array}{c}
 \|W\|_{\op}\le L/2,\quad
 \lambda_j^+(W)\ge c_0((k+j)/d)^2\ \text{for all }j,\\[1mm]
 A_K/d^2\le\mu:=\lambda_{\min}^+(W)\le A_K^+/d^2,\\[1mm]
 \nu(P):=\lambda_{\min}(QWQ|_{\range(Q)})\le B_K/d^2
 \end{array}
 \right]
 \ge\frac13+\rho_K.
 \label{eq:finite-k-event}
\end{equation}
The constants can be chosen so that $B_K<(1-\tau)A_K$ for a universal
$\tau>0$.  On the same event, if
$s_1$ is the least positive eigenvalue of $S$ and
$\theta_1$ is the least eigenvalue of
$QSQ|_{\range(Q)}$, then
\begin{equation}
 s_1\ge\frac{\mathsf A_K}{d^2},
 \qquad
 \theta_1\le\frac{\mathsf B_K}{d^2},
 \qquad
 0<\mathsf B_K<\mathsf A_K,
 \label{eq:S-one-mode-separation}
\end{equation}
for universal constants depending only on the fixed $K$.
\end{lemma}

\begin{proof}
We first isolate the two hard-edge events.  Write $n=d-k$ and couple the
standard Gaussian matrix $G_{d,n}$ with its first $n+1$ rows
$G_{n+1,n}$.  Adding rows can only increase the least singular value, and
$\sqrt d\ge\sqrt{n+1}$, so
\begin{equation}
 \sqrt d\,s_{\min}(G_{d,n})
 \ge
 \sqrt{n+1}\,s_{\min}(G_{n+1,n}).
 \label{eq:row-monotonicity}
\end{equation}
Since $X=d^{-1/2}G_{d,n}$,
\[
 d^2\mu
 =\bigl(\sqrt d\,s_{\min}(G_{d,n})\bigr)^2.
\]
Consequently, uniformly for $1\le k<K$ and all sufficiently large $d$,
\begin{equation}
 \Prb[d^2\mu<a_0^2]
 \le F_1(a_0)+o_K(1).
 \label{eq:global-edge-prob}
\end{equation}

Let $g:=F_0(b_0)-F_1(a_0)-1/3>0.036$ and fix $\delta_*=g/20$.
The required upper tail has a nonasymptotic proof. Append $k$ independent
Gaussian columns to $G_{d,d-k}$ to form a square matrix $J_d$.
Principal-submatrix interlacing and \citet[Theorem~1.6]{Wei17} give
\[
 s_{\min}(G_{d,d-k})\le s_{k+1}^{\up}(J_d),\qquad
 \Prb[d^2\mu>C_1^2t^2(k+1)^2]\le e^{-C_2t(k+1)}.
\]
Choose a fixed $t>1$ with $e^{-2C_2t}\le\delta_*$ and then
$R_K>\max\{a_0^2,C_1^2t^2K^2\}$. Thus
\begin{equation}
 \Prb[d^2\mu>R_K]\le\delta_*
 \label{eq:global-edge-upper}
\end{equation}
for all such $k$ and $d$. Increase $d_K$ so that the error in \eqref{eq:global-edge-prob} is at most $\delta_*$. This step requires no limiting law at rectangularity $k>1$.

For almost every actual transcript, the construction in \Cref{thm:intermediate} couples the compression with a fresh square $h\times h$ Gaussian matrix, where $h=d-k-q$, and gives
\[
 \nu(P)\le d^{-1}s_{\min}(G_h\Sigma)^2
 \le d^{-1}s_{\min}(G_h)^2,\qquad \|\Sigma\|_{\op}\le1.
\]
Consequently,
\[
 \Prb[\nu(P)\le b_0^2/d^2\mid\cF_q]
 \ge\Prb[\sqrt h\,s_{\min}(G_h)\le b_0\sqrt{h/d}].
\]
Take $d_K\ge d_0$ and $M_K\ge4$, so that \Cref{lem:rect-spectrum} applies.
Choose a fixed $\zeta>0$ sufficiently small that
$F_0(b_0\sqrt{1-\zeta})\ge F_0(b_0)-\delta_*$. Choose $\alpha_K$ small and $M_K$ large so that $\alpha_K+1/M_K\le\zeta$. Then $h/d\ge1-\zeta$. Square hard-edge convergence, with $d_K$ increased if necessary, makes the last probability at least $F_0(b_0)-2\delta_*$ uniformly in $h$ in this range and almost surely in the transcript. This is a fixed continuity loss plus a convergence error; the continuity loss is not claimed to vanish with $d$.

After integrating the conditional bound, a union bound gives
\[
 \Prb[d^2\mu\ge a_0^2,\ d^2\nu(P)\le b_0^2]
 \ge F_0(b_0)-F_1(a_0)-3\delta_*.
\]
Intersect with \eqref{eq:global-edge-upper} and the joint spectral event in \Cref{lem:rect-spectrum}. The total probability remains at least
\[
 \frac13+g-4\delta_*-0.005
 >\frac13+\frac g2.
\]
Indeed, $g-4\delta_*-0.005=0.8g-0.005>g/2$ because $g>0.036$. No independence of these events is used. Thus one may take $\rho_K=g/2$, $A_K=a_0^2$, $B_K=b_0^2$, and $A_K^+=R_K$. Their ratio is $B_K/A_K=225/256$, so any fixed $\tau<31/256$ satisfies the required gap.

Finally, on $\mu\le A_K^+/d^2$,
\[
 s_1
 =\frac{2\mu}{L}-\frac{\mu^2}{L^2}
 \ge\frac{2(1-o(1))A_K}{Ld^2},
\]
whereas
\[
 \theta_1
 \le\frac2L\nu(P)
 \le\frac{2B_K}{Ld^2}.
\]
Because $A_K>B_K$, increasing $d_K$ yields constants
$\mathsf A_K>\mathsf B_K>0$ satisfying
\eqref{eq:S-one-mode-separation}.
\end{proof}

\begin{remark}
The probability threshold in \Cref{lem:finite-k-separation} is the
reason the explicit hard-edge laws matter.  Generic ``constant probability''
upper and lower bounds would not suffice for a $2/3$-success minimax theorem.
The strict CDF gap in \eqref{eq:CDF-gap} leaves enough probability after all
auxiliary spectral events are intersected.
\end{remark}

\section{From Spectral Profiles to Schatten Error}\label{sec:heat}

Throughout this section, $2\le p<\infty$ and we work on the spectral event $\|W\|_{\op}\le L/2$. Thus $0\preceq S\preceq(3/4)I$. Let $P$ be the output projector, $Q=I-P$, and let
$0\le\theta_1\le\cdots\le\theta_r$ be the eigenvalues of
$QSQ|_{\range(Q)}$, where $r=d-k$.  Because
$QH^2Q=I_{\range(Q)}-QSQ$,
\begin{equation}
 R_p(P):=\|HQ\|_{S_p}^p
 =\sum_{j=1}^{r}(1-\theta_j)^{p/2}.
 \label{eq:R}
\end{equation}
The optimal rank-$k$ projector removes the $k$ unit eigenvalues of $H$, so
\begin{equation}
 Z_p:=\OPT_{p,k}(H)^p
 =\sum_{j=1}^{r}(1-s_j)^{p/2}.
 \label{eq:Z}
\end{equation}
Poincare separation for the $r$-dimensional compression of the
$d$-dimensional matrix $S$, whose first $k$ eigenvalues are zero, gives
\begin{equation}
 \theta_j\le s_j
 \qquad(1\le j\le r).
 \label{eq:theta-s}
\end{equation}
Thus every summand in $R_p(P)-Z_p$ is nonnegative.

\begin{lemma}[Optimal heat-trace bounds]
\label[lemma]{lem:heat-upper}
On the event of \Cref{lem:rect-spectrum}, for every $p\ge2$,
\begin{align}
 Z_p
 &\le C\frac{d}{\sqrt p},
 \label{eq:heat1}\\
 Z_p
 &\le
 C\left(1+\frac{d^2}{pk}\right)
 \exp\!\left(-c\frac{pk^2}{d^2}\right).
 \label{eq:heat2}
\end{align}
\end{lemma}

\begin{proof}
Since $(1-x)^{p/2}\le e^{-px/2}$ and
$s_j\ge c(k+j)^2/d^2$,
\[
 Z_p\le\sum_{j\ge1}
 \exp\!\left[-c\frac{p(k+j)^2}{d^2}\right].
\]
Bounding the sum by a Gaussian integral gives \eqref{eq:heat1}.  Starting the
same integral at $k$ and using the standard Gaussian-tail estimate gives
\eqref{eq:heat2}.
\end{proof}

Set the dimensionless hard-edge temperature
\begin{equation}
 a^2:=\frac{pk^2}{d^2},
 \qquad D:=\frac{d}{\sqrt p}=\frac{k}{a}.
 \label{eq:aD}
\end{equation}

\begin{samepage}
\begin{lemma}[Large-rank multimode comparison]
\label[lemma]{lem:comparison-large}
Assume the lower-spectrum event of \Cref{lem:rect-spectrum} and the fresh
intermediate event of \Cref{thm:intermediate}, with
$k+q\le d/2$.  There exist universal constants
$K_\star,A,c,C>0$ such that whenever $k\ge K_\star$:
\begin{enumerate}[label=(\roman*),leftmargin=2em]
\item If $a^2\le A$, then
\begin{equation}
 R_p(P)-Z_p\ge c\frac{pk^3}{d^2}.
 \label{eq:linear-excess}
\end{equation}
\item If $a^2>A$, then
\begin{equation}
 \frac{R_p(P)}{Z_p}
 \ge\exp(ca^2).
\label{eq:exponential-ratio}
\end{equation}
\end{enumerate}
\end{lemma}
\end{samepage}

\begin{proof}
For part (i), take $\ell=\lfloor\gamma k\rfloor$ with a sufficiently small
universal $\gamma>0$.  Once $k\ge K_\star$, we have
$\ell\ge\gamma k/2$; also $\ell\le k\le h$ after choosing $\gamma\le1$.  For $j\le\ell$, the fresh bound and the rectangular
lower bound give
\[
 \theta_j\le C\gamma^2\frac{k^2}{d^2},
 \qquad
 s_j\ge c_1\frac{k^2}{d^2}.
\]
Choose $\gamma$ so that the first coefficient is at most $c_1/2$.  For
$f_p(x)=(1-x)^{p/2}$, integrate $-f_p'$ only from $\theta_j$ to
$c_1k^2/d^2$, a subinterval of $[\theta_j,s_j]$.  Since
$pk^2/d^2=a^2\le A$, the derivative is bounded below by a constant multiple
of $p$ on this interval.  Thus
\[
 f_p(\theta_j)-f_p(s_j)
 \ge c_Ap\frac{k^2}{d^2}.
\]
Summing over $\ell=\Theta(k)$ indices proves
\eqref{eq:linear-excess}.

For part (ii), first suppose $D=d/\sqrt p$ exceeds a sufficiently large
constant. Put $J=\lfloor\beta D\rfloor$ with $\beta>0$ small. Since $p\ge2$ and $h\ge d/2$, choosing $\beta\le1/2$ ensures $J\le h$. Raise the threshold on $D$ so that $J\ge\beta D/2$.  For
$j\le J$, the fresh bound gives $p\theta_j\le C\beta^2$, and hence
$(1-\theta_j)^{p/2}\ge c$.  Therefore
\[
 R_p(P)\ge cD.
\]
\Cref{lem:heat-upper} gives
\[
 Z_p\le C(1+D/a)e^{-ca^2}.
\]
It follows that
\[
 \frac{R_p(P)}{Z_p}
 \ge c\frac{D}{1+D/a}e^{ca^2}
 \ge c\min\{D,a\}e^{ca^2}.
\]
Taking $A$ and the lower threshold on $D$ large absorbs the prefactor.

If $D$ is bounded, use the first fresh eigenvalue:
\[
 R_p(P)
 \ge(1-C/d^2)^{p/2}
 \ge c\exp(-Cp/d^2)
 =c\exp(-Ca^2/k^2).
\]
On the other hand, \eqref{eq:heat2} gives
$Z_p\le C\exp(-ca^2)$.  Choose $K_\star$ so large that
$c-C/K_\star^2>0$, and then enlarge $A$ to absorb the remaining constants.
This proves \eqref{eq:exponential-ratio}.
\end{proof}

\begin{samepage}
\begin{lemma}[Bounded-rank one-mode comparison]
\label[lemma]{lem:comparison-small}
Fix the universal $K_\star$ from
\Cref{lem:comparison-large}.  On the event of
\Cref{lem:finite-k-separation} with $K=K_\star$, assume
$1\le k<K_\star$ and put
\[
 u:=\frac{p}{d^2}.
\]
There exist universal constants $U_\star,c>0$ such that:
\begin{enumerate}[label=(\roman*),leftmargin=2em]
\item If $u\le U_\star$, then
\begin{equation}
 R_p(P)-Z_p\ge c\frac{p}{d^2}.
 \label{eq:small-k-linear}
\end{equation}
\item If $u>U_\star$, then
\begin{equation}
 \frac{R_p(P)}{Z_p}\ge\exp(cu).
\label{eq:small-k-exp}
\end{equation}
\end{enumerate}
\end{lemma}
\end{samepage}

\begin{proof}
Write
\[
 \mathsf A:=\mathsf A_{K_\star},
 \qquad
 \mathsf B:=\mathsf B_{K_\star},
 \qquad
 \Delta:=\mathsf A-\mathsf B>0.
\]
On the separation event,
\[
 \theta_1\le\mathsf B/d^2,
 \qquad
 s_1\ge\mathsf A/d^2.
\]
Since every summand in $R_p(P)-Z_p$ is nonnegative,
\[
 R_p(P)-Z_p
 \ge
 f_p(\mathsf B/d^2)-f_p(\mathsf A/d^2),
 \qquad f_p(x)=(1-x)^{p/2}.
\]
If $u\le U_\star$, the mean-value theorem and
$p\mathsf A/d^2\le \mathsf A U_\star$ give
\[
 f_p(\mathsf B/d^2)-f_p(\mathsf A/d^2)
 \ge c_{U_\star}\frac{p\Delta}{d^2},
\]
which proves \eqref{eq:small-k-linear}.

Now suppose $u>U_\star$.  Increase the universal dimension threshold in
\Cref{lem:finite-k-separation} so that
\[
 (1-\mathsf B/d^2)^{p/2}
 \ge
 \exp\!\left[
   -\frac{\mathsf B+\Delta/4}{2}\,u
 \right].
\]
Thus
\begin{equation}
 R_p(P)
 \ge
 \exp\!\left[
   -\frac{\mathsf B+\Delta/4}{2}\,u
 \right].
 \label{eq:R-small-k-lower}
\end{equation}
Choose a universal integer $J$ so large that the profile constant $c_1$ in
\eqref{eq:S-spectrum} satisfies
\[
 \frac{c_1(J+1)^2}{4}
 \ge
 \frac{\mathsf B+\Delta/4}{2}+\frac{\Delta}{4}.
\]
For the first $J$ optimal modes, $s_j\ge s_1\ge\mathsf A/d^2$, while for
$j>J$, $s_j\ge c_1j^2/d^2$.  Hence
\begin{align*}
 Z_p
 &\le
 J e^{-\mathsf A u/2}
 +\sum_{j>J}e^{-c_1u j^2/2}\\
 &\le
 C_J\exp\!\left[
   -\left(
      \frac{\mathsf B+\Delta/4}{2}+\frac{\Delta}{4}
    \right)u
 \right]
\end{align*}
for all sufficiently large $U_\star$.  Combining this with
\eqref{eq:R-small-k-lower} gives
$R_p(P)/Z_p\ge C_J^{-1}e^{\Delta u/4}$.  Enlarging $U_\star$ once more
absorbs $C_J$ and proves \eqref{eq:small-k-exp}.
\end{proof}

The linear comparison controls the total contribution of the first modes. The exponential comparison explains the spectral behavior at larger $p$: an output complement retaining an eigenvalue below the optimal edge has a multiplicatively larger residual power.

\section{Completing the Uniform Lower Bounds}\label{sec:assembly}

\begin{proof}[Proof of the finite-$p$ lower bound in
\Cref{thm:main}]
By Yao's principle, fix a deterministic adaptive algorithm and draw the hard
instance.  Let $K_\star$ be the constant in \Cref{lem:comparison-large}, enlarged also to satisfy the large-rank spectral comparison below. Fix this value before choosing any bounded-rank constants, and define
\[
 R_*:=\min\{N,k\Phi(p,\eps)\}.
\]
Choose the spectral and comparison constants first. Next choose $\eta$ small enough for the inequalities below, then $M_0$ large enough for the hard-dimension thresholds, and finally the query fraction $\alpha$ small enough for the posterior events. All these constants are universal.

First suppose
\begin{equation}
 R_*\le M_0k
 \label{eq:baseline-branch}
\end{equation}
where the universal $M_0$ is chosen after the dimension constant $\eta$.  Apply
\Cref{lem:range-baseline} in square dimension $N$ and, when the requested
matrix is rectangular, embed it by zero padding.  This gives the exact lower
bound
\[
 q\ge k\ge R_*/M_0,
\]
so $q=\Omega(R_*)$ with no restriction on the ratio $k/N$.  This branch covers exactly the case $R_*=O(k)$, including near-full rank. When $R_*=N\gg k$, the following Wishart branch is still necessary.

Assume henceforth that $R_*>M_0k$.  Choose a sufficiently small universal
$\eta>0$ and set
\begin{equation}
 d:=\lfloor\eta R_*\rfloor.
 \label{eq:dchoice}
\end{equation}
Taking $M_0$ sufficiently large after $\eta$ is fixed ensures $d\ge\eta R_*/2$ and simultaneously that
\[
 d\ge M_{K_\star}k,\qquad d\ge d_{K_\star},\qquad d\ge d_0,
 \qquad k\le d/4.
\]
Also $d\le N$, so the $d\times d$ hard matrix can be embedded in the requested
dimensions by \Cref{lem:padding}.  Suppose toward a contradiction that
the algorithm uses
\[
 q\le\alpha d,
\]
where $\alpha>0$ is chosen below all constants required by
\Cref{thm:intermediate} and
\Cref{lem:finite-k-separation}.  Then $k+q\le d/2$.

The definition of $\Phi$ gives the two scale inequalities
\begin{align}
 d
 &\le \eta k\Phi(p,\eps)
 \le \eta\frac{k}{\sqrt\eps},
 \label{eq:d-spectral-cap}\\
 d^3
 &\le \eta^3k^3\Phi(p,\eps)^3
 \le \eta^3\frac{k^3\sqrt p}{\eps}.
 \label{eq:d-cubic}
\end{align}
We split according to the target rank.

\medskip
\noindent\emph{Case 1: $k\ge K_\star$.}
On the intersection of the uniform rectangular-spectrum event and the
conditional intermediate hard-edge event, which has probability at least
$0.9$, \Cref{lem:comparison-large} applies.  Set
\[
 a^2:=\frac{pk^2}{d^2}.
\]
Equation~\eqref{eq:d-spectral-cap} implies
\begin{equation}
 p\eps\le\eta^2a^2.
 \label{eq:peps-a}
\end{equation}
If $a^2\le A$, then $p\eps\le\eta^2A$.  For $\eta$ small,
$e^{p\eps}-1\le C p\eps$, and relative-error success would imply
\[
 R_p(P)-Z_p
 \le\bigl((1+\eps)^p-1\bigr)Z_p
 \le C p\eps\,\frac{d}{\sqrt p}
 =C\eps d\sqrt p.
\]
\Cref{lem:comparison-large}(i) gives the opposite bound
$c pk^3/d^2$.  By \eqref{eq:d-cubic}, the ratio of the lower bound to the
upper bound is at least $c/(C\eta^3)$, a contradiction for sufficiently small
$\eta$.

If $a^2>A$, success would give
\[
 \frac{R_p(P)}{Z_p}
 \le(1+\eps)^p
 \le e^{p\eps}
 \le e^{\eta^2a^2},
\]
whereas \Cref{lem:comparison-large}(ii) gives
$R_p(P)/Z_p\ge e^{ca^2}$.  Taking $\eta^2<c$ again gives a contradiction.
Thus every $q\le\alpha d$ algorithm fails on an event of probability at least
$0.9$ in this case.

\medskip
\noindent\emph{Case 2: $1\le k<K_\star$.}
Use the finite-codimension event from
\Cref{lem:finite-k-separation} with $K=K_\star$.  Its probability is at
least $1/3+\rho_{K_\star}$, uniformly over the choice of adaptive algorithm; the joint event is unconditional.  Put
\[
 u:=\frac{p}{d^2}.
\]
Equation~\eqref{eq:d-spectral-cap} now gives
\begin{equation}
 p\eps\le\eta^2k^2u\le\eta^2K_\star^2u.
 \label{eq:peps-u}
\end{equation}
If $u\le U_\star$, choose $\eta$ so that
$\eta^2K_\star^2U_\star$ is a sufficiently small constant.  Then relative
success and \Cref{lem:heat-upper} imply
\[
 R_p(P)-Z_p
 \le C\eps d\sqrt p.
\]
On the other hand, \Cref{lem:comparison-small}(i) gives
$R_p(P)-Z_p\ge cp/d^2$.  Their ratio satisfies
\[
 \frac{cp/d^2}{C\eps d\sqrt p}
 =\frac{c\sqrt p}{C\eps d^3}
 \ge\frac{c}{C\eta^3k^3}
 \ge\frac{c}{C\eta^3K_\star^3},
\]
where we used \eqref{eq:d-cubic}.  This is a contradiction after decreasing
$\eta$.

If $u>U_\star$, success and \eqref{eq:peps-u} give
\[
 \frac{R_p(P)}{Z_p}
 \le e^{p\eps}
 \le e^{\eta^2K_\star^2u},
\]
while \Cref{lem:comparison-small}(ii) gives
$R_p(P)/Z_p\ge e^{cu}$.  Choose
$\eta^2K_\star^2<c$.  Hence every $q\le\alpha d$ algorithm fails on an event
of probability strictly larger than $1/3$ also in the bounded-rank case.

Combining both cases, no deterministic $q\le\alpha d$ algorithm has success
probability $2/3$ under the relevant hard distribution.  Yao's minimax
principle therefore gives
\[
 q=\Omega(d)=\Omega(R_*).
\]
The hard matrix is symmetric, so the reduction is valid for the full
two-sided oracle.
\end{proof}

\begin{proof}[Proof of the spectral endpoint in
\Cref{thm:main}]
Define
\[
 R_\infty:=\min\left\{N,\frac{k}{\sqrt\eps}\right\}.
\]
The branch $R_\infty\le M_0k$ follows from
\Cref{lem:range-baseline}, exactly as in
\eqref{eq:baseline-branch}, and does not require $k\le cN$.  Otherwise take
\[
 d=\lfloor\eta R_\infty\rfloor
\]
with the same hierarchy of sufficiently small and large universal constants.
Then
\begin{equation}
 d\le\eta\frac{k}{\sqrt\eps},
 \qquad
 \eps\le\eta^2\frac{k^2}{d^2}.
 \label{eq:spectral-scale}
\end{equation}
Consider any deterministic algorithm with $q\le\alpha d$ queries.

If $k\ge K_\star$, \Cref{lem:rect-spectrum} and the $j=1$ case of
\Cref{thm:intermediate} give, on an event of probability at least
$0.9$,
\[
 \mu\ge c_0\frac{k^2}{d^2},
 \qquad
 \nu(P)\le\frac{C}{d^2}.
\]
The common choice of $K_\star$, made before fixing the bounded-rank constants, ensures $C<(9/16)c_0K_\star^2$.  By \eqref{eq:spectral-scale}, decreasing $\eta$
ensures $\eps\le(\mu/L)/8$.  Spectral relative-error success would then imply,
by \eqref{eq:spectral-certificate-coarse},
\[
 \nu(P)\ge\frac{9}{16}\mu
 \ge\frac{9c_0}{16}\frac{k^2}{d^2},
\]
contradicting the fresh hard-edge upper bound.

It remains to consider $1\le k<K_\star$.  On the event in
\Cref{lem:finite-k-separation},
\[
 A_{K_\star}/d^2\le\mu\le A_{K_\star}^+/d^2,
 \qquad
 \nu(P)\le B_{K_\star}/d^2,
 \qquad
 B_{K_\star}<(1-\tau)A_{K_\star}.
\]
Fix $\eta_0<\tau/2$.  For $d$ above the lemma's dimension threshold,
$\delta=\mu/L\le\eta_0$; and by \eqref{eq:spectral-scale}, choosing the hard
dimension constant $\eta$ sufficiently small gives
$\eps\le\eta_0\delta$, uniformly over $k<K_\star$.  The sharp conclusion
\eqref{eq:spectral-certificate} would force
\[
 \nu(P)\ge(1-2\eta_0)\mu
 >\frac{B_{K_\star}}{d^2},
\]
again a contradiction.  The finite-codimension event has probability
strictly greater than $1/3$.

Yao's principle now yields
\[
 q=\Omega(R_\infty).
\]
Together with spectral Krylov iteration and exact recovery, this proves
\eqref{eq:spectral-main}.
\end{proof}

\section{Inherited Upper Bounds and the Transition}
\label{sec:upper}
For completeness, we specify how the known algorithms cover the same oracle, output, and probability convention. Exact recovery takes $N$ products by querying the standard basis on the smaller side. All subsequent singular-value computations are uncharged. For fixed finite $p$, the Schatten approximation algorithm of \citet{BCW22} gives the upper bounds in \Cref{prop:small-p} and \Cref{cor:fixed} with success probability exceeding $2/3$.

Dimension logarithms can use $N$ by applying the algorithm with the smaller dimension on the right. If $n>m$, run it on $A^\top$ to obtain a rank-$k$ left projector $UU^\top$ for $A$. With $k$ further products, form $A^\top U$ and let $P$ project onto its range, completed to rank $k$ when necessary. Since $U^\top A(I-P)=0$,
\[
 \|A(I-P)\|_{S_p}
 =\|(I-UU^\top)A(I-P)\|_{S_p}
 \le\|(I-UU^\top)A\|_{S_p}.
\]
The left and right optimal residuals agree. This conversion preserves the guarantee and adds only $k$ queries.

For the uniform result, a useful explicit form of the BCW guarantee, including this harmless conversion cost, is
\begin{equation}
 O\!\left(k\log(N/\eps)
       \bigl(p^{1/6}\eps^{-1/3}+\sqrt p\bigr)\right)
 \label{eq:bcw-explicit}
\end{equation}
products for $1\le p\lesssim\log(N)/\eps$, with success probability at least $9/10$ \citep[Theorem~5.1 of the arXiv full version; Theorem~4.2 of the proceedings version]{BCW22}. The second term must be accounted for when $p$ grows. If $p\eps\le1$, then $\sqrt p\le p^{1/6}\eps^{-1/3}$, so \eqref{eq:bcw-explicit} has polynomial rate $k\Phi$. If $1\le p\eps\lesssim\log N$, both terms are bounded by $k\eps^{-1/2}$ times logarithmic factors: their additional factors are respectively $(p\eps)^{1/6}$ and $(p\eps)^{1/2}$. Hence this corridor also has rate $\wtO(k\Phi)$.

If $p\ge C\log N/\eps$ for a sufficiently large constant $C$, then every residual matrix $B$ has rank at most $N$ and
\[
 \|B\|_{S_p}\le N^{1/p}\|B\|_{\op}.
\]
Run spectral block Krylov iteration at accuracy $\eps/4$. Choose $C$ so that $N^{1/p}\le1+\eps/4$ for $\eps\le1$. Its output satisfies
\[
 \|A(I-P)\|_{S_p}
 \le(1+\eps/4)^2\OPT_{\infty,k}(A)
 \le(1+\eps)\OPT_{p,k}(A).
\]
The required number of products is $\wtO(k/\sqrt\eps)$ \citep{MM15}. The same algorithm directly treats $p=\infty$. These algorithms have predetermined query budgets and constant success probability at least $2/3$. Selecting the smaller of the applicable budget and $N$ proves the upper bounds in \Cref{thm:main,prop:small-p}.

\section{Discussion}
\label{sec:discussion}
The two proofs isolate complementary information barriers. Posterior overlap explains the multiplicative rank factor with a short argument valid for every fixed finite Schatten parameter. Spectral-profile persistence identifies the additional cost of increasing $p$ and its transition to the spectral rate. The upper bounds are inherited; the contribution is the adaptive lower bound with joint parameter dependence.

For $k=N-r$, the theorem gives linear complexity up to logarithms whenever the target rate reaches $N$. If $r\ll N$, then $k=\Theta(N)$ and the exact range-recovery baseline suffices. If saturation instead comes from stringent accuracy with $k\ll N$, the rectangular Wishart construction supplies the linear lower bound. In both cases, the worst-case multiplicative model is essential: exactly rank-$k$ inputs have zero optimal residual. A promise of full rank or a positive lower bound on the optimal tail would require a separate analysis of the baseline branch.

The results leave exact logarithmic factors open and concern exact real-valued products with uncharged computation. They do not characterize noisy or finite-precision oracles. \Cref{prop:small-p} covers each fixed $1\le p<2$, and \Cref{cor:fixed} summarizes the consequence for all fixed finite norms. We do not claim constants uniform over $1\le p<2$. The hard-edge probability argument is specifically calibrated to success probability $2/3$; no confidence-dependent minimax characterization is asserted.

\clearpage
\bibliographystyle{unsrtnat}
\setlength{\bibsep}{3pt}
\bibliography{references}

@inproceedings{BCW22,
  author    = {Ainesh Bakshi and Kenneth L. Clarkson and David P. Woodruff},
  title     = {Low-Rank Approximation with {$1/\epsilon^{1/3}$} Matrix-Vector Products},
  booktitle = {Proceedings of the 54th Annual ACM SIGACT Symposium on Theory of Computing (STOC)},
  pages     = {1130--1143},
  year      = {2022},
  doi       = {10.1145/3519935.3519988},
  url       = {https://arxiv.org/abs/2202.05120}
}

@inproceedings{BN23,
  author    = {Ainesh Bakshi and Shyam Narayanan},
  title     = {Krylov Methods are (Nearly) Optimal for Low-Rank Approximation},
  booktitle = {Proceedings of the 64th IEEE Symposium on Foundations of Computer Science (FOCS)},
  pages     = {2093--2101},
  year      = {2023},
  doi       = {10.1109/FOCS57990.2023.00128},
  url       = {https://arxiv.org/abs/2304.03191}
}

@inproceedings{AmselEtAl26,
  author    = {Noah Amsel and Pratyush Avi and Tyler Chen and Feyza Duman Keles and Chinmay Hegde and Christopher Musco and Cameron Musco and David Persson},
  title     = {Query Efficient Structured Matrix Learning},
  booktitle = {Proceedings of the 39th Annual Conference on Learning Theory (COLT)},
  series    = {Proceedings of Machine Learning Research},
  volume    = {336},
  pages     = {158--194},
  year      = {2026},
  url       = {https://proceedings.mlr.press/v336/amsel26a.html}
}

@inproceedings{BHSW20,
  author    = {Mark Braverman and Elad Hazan and Max Simchowitz and Blake Woodworth},
  title     = {The Gradient Complexity of Linear Regression},
  booktitle = {Proceedings of the 33rd Annual Conference on Learning Theory (COLT)},
  series    = {Proceedings of Machine Learning Research},
  volume    = {125},
  pages     = {627--647},
  year      = {2020},
  url       = {https://proceedings.mlr.press/v125/braverman20a.html}
}

@article{RV09,
  author  = {Mark Rudelson and Roman Vershynin},
  title   = {The Smallest Singular Value of a Random Rectangular Matrix},
  journal = {Communications on Pure and Applied Mathematics},
  volume  = {62},
  number  = {12},
  pages   = {1707--1739},
  year    = {2009},
  doi     = {10.1002/cpa.20294},
  url     = {https://arxiv.org/abs/0802.3956}
}

@inproceedings{MM15,
  author    = {Cameron Musco and Christopher Musco},
  title     = {Randomized Block Krylov Methods for Stronger and Faster Approximate Singular Value Decomposition},
  booktitle = {Advances in Neural Information Processing Systems 28 (NeurIPS)},
  pages     = {1396--1404},
  year      = {2015},
  url       = {https://proceedings.neurips.cc/paper/2015/hash/1efa39bcaec6f3900149160693694536-Abstract.html}
}

@inproceedings{CEMMR26,
  author    = {Tyler Chen and Ethan N. Epperly and Raphael A. Meyer and Christopher Musco and Akash Rao},
  title     = {Does Block Size Matter in Randomized Block Krylov Low-Rank Approximation?},
  booktitle = {Proceedings of the 37th Annual ACM-SIAM Symposium on Discrete Algorithms (SODA)},
  pages     = {1026--1046},
  year      = {2026},
  doi       = {10.1137/1.9781611978971.42},
  url       = {https://arxiv.org/abs/2508.06486}
}

@article{Edelman88,
  author  = {Alan Edelman},
  title   = {Eigenvalues and Condition Numbers of Random Matrices},
  journal = {SIAM Journal on Matrix Analysis and Applications},
  volume  = {9},
  number  = {4},
  pages   = {543--560},
  year    = {1988},
  doi     = {10.1137/0609045}
}

@article{Edelman91,
  author  = {Alan Edelman},
  title   = {The Distribution and Moments of the Smallest Eigenvalue of a Random Matrix of Wishart Type},
  journal = {Linear Algebra and its Applications},
  volume  = {159},
  pages   = {55--80},
  year    = {1991},
  doi     = {10.1016/0024-3795(91)90076-9}
}

@article{NF98,
  author  = {Taro Nagao and Peter J. Forrester},
  title   = {The Smallest Eigenvalue Distribution at the Spectrum Edge of Random Matrices},
  journal = {Nuclear Physics B},
  volume  = {509},
  number  = {3},
  pages   = {561--598},
  year    = {1998},
  doi     = {10.1016/S0550-3213(97)00670-6}
}

@article{Wei17,
  author  = {Feng Wei},
  title   = {Upper Bound for Intermediate Singular Values of Random Matrices},
  journal = {Journal of Mathematical Analysis and Applications},
  volume  = {445},
  number  = {2},
  pages   = {1530--1547},
  year    = {2017},
  doi     = {10.1016/j.jmaa.2016.08.007},
  url     = {https://arxiv.org/abs/1606.03931}
}

@article{HMT11,
  author = {Nathan Halko and Per-Gunnar Martinsson and Joel A. Tropp},
  title = {Finding Structure with Randomness: Probabilistic Algorithms for Constructing Approximate Matrix Decompositions},
  journal = {SIAM Review},
  volume = {53},
  number = {2},
  pages = {217--288},
  year = {2011},
  doi = {10.1137/090771806},
  url = {https://arxiv.org/abs/0909.4061}
}

@inproceedings{MMM24,
  author = {Raphael A. Meyer and Cameron Musco and Christopher Musco},
  title = {On the Unreasonable Effectiveness of Single Vector Krylov Methods for Low-Rank Approximation},
  booktitle = {Proceedings of the Annual ACM-SIAM Symposium on Discrete Algorithms (SODA)},
  pages = {811--845},
  year = {2024},
  doi = {10.1137/1.9781611977912.32},
  url = {https://arxiv.org/abs/2305.02535}
}

@inproceedings{KW24,
  author = {Praneeth Kacham and David P. Woodruff},
  title = {Faster Algorithms for {Schatten-$p$} Low Rank Approximation},
  booktitle = {Approximation, Randomization, and Combinatorial Optimization. Algorithms and Techniques (APPROX/RANDOM)},
  series = {Leibniz International Proceedings in Informatics},
  volume = {317},
  pages = {55:1--55:19},
  year = {2024},
  doi = {10.4230/LIPIcs.APPROX/RANDOM.2024.55},
  url = {https://arxiv.org/abs/2407.11959}
}

@inproceedings{KW23,
  author = {Praneeth Kacham and David P. Woodruff},
  title = {Lower Bounds on Adaptive Sensing for Matrix Recovery},
  booktitle = {Advances in Neural Information Processing Systems},
  volume = {36},
  year = {2023},
  url = {https://arxiv.org/abs/2311.17281}
}
\end{document}